\documentclass[a4paper,pdftex,reqno,10pt]{amsart}

\allowdisplaybreaks[1]

\usepackage{geometry}
\usepackage{graphicx}
\usepackage{enumerate}
\usepackage{courier}
\usepackage{times}
\usepackage{amssymb}
\usepackage{amsmath}
\usepackage{footmisc}

\usepackage{hyperref}

\theoremstyle{plain}
\newtheorem{Theorem}{Theorem}
\newtheorem{Proposition}[Theorem]{Proposition}

\newtheorem{Lemma}[Theorem]{Lemma}

\newenvironment{Proof}
{\begin{trivlist}\item[]{{\sc Proof.}}}{\hfill{$\square$}\noindent\end{trivlist}}

\newcommand{\gaussm}[3]{\genfrac{[}{]}{0pt}{}{#1}{#2}_{#3}}

\newcommand{\F}[2]{\mathbb{F}_{#2}^{#1}}

\newcommand{\SV}[2]{\mathbb{F}_{#2}^{#1}}

\theoremstyle{definition}
\newtheorem{Definition}[Theorem]{Definition}

\newtheorem{Example}[Theorem]{Example}

\theoremstyle{remark}

\begin{document}


\title{Heden's bound on the tail of a vector space partition}

\author[Sascha Kurz]{Sascha Kurz$^\star$}
\address{Department of Mathematics, University of Bayreuth, 95440 Bayreuth, Germany}
\email{sascha.kurz@uni-bayreuth.de}
\thanks{$^\star$ 
Grant KU 2430/3-1 -- 
\emph{Integer Linear Programming Models for Subspace Codes and Finite Geometry} 
-- German Research Foundation.}

\date{}

\begin{abstract}
  A vector space partition of $\F{v}{q}$ is a collection of subspaces such that every non-zero 
  vector is contained in a unique element. We 
  improve a lower bound of Heden, in a subcase, on the 
  number of elements of the smallest occurring dimension in a vector space partition. To this end, we introduce the notion 
  of $q^r$-divisible sets of $k$-subspaces in $\F{v}{q}$. By geometric arguments we obtain non-existence results for 
  these objects, which then imply the improved result of Heden. 
\end{abstract}

\maketitle

\section{Introduction}
\noindent
Let $q>1$ be a prime power, $\mathbb{F}_q$ be the finite filed with $q$ elements, and $v$ a positive integer. A \emph{vector space partition} 
$\mathcal{P}$ of $\F{v}{q}$ is a collection of subspaces with the property that every non-zero vector is contained in a unique member of 
$\mathcal{P}$. If $\mathcal{P}$  contains $m_d$ subspaces of dimension $d$, then $\mathcal{P}$ is of type $k^{m_k}\ldots 1^{m_1}$. We may leave 
out some of the cases with $m_d=0$. Subspaces of dimension~$d$ are also called \emph{$d$-subspaces}. $1$-subspaces are called \emph{points}, 
$(v-1)$-subspaces are called \emph{hyperplanes}, and 
each $k$-subspace contains $\gaussm{k}{1}{q}:=\frac{q^k-1}{q-1}$ points. So, in a vector space partition $\mathcal{P}$ each point of the ambient 
space $\F{v}{q}$ is covered by exactly one point of one of the elements of $\mathcal{P}$. An example of a vector space partition is given 
by a \emph{$k$-spread} in $\F{v}{q}$, where $\gaussm{v}{1}{q}/\gaussm{k}{1}{q}$ $k$-subspaces partition the set of points of $\F{v}{q}$. The 
corresponding type is given by $k^{m_k}$, where $m_k=\gaussm{v}{1}{q}/\gaussm{k}{1}{q}$. If $k$ divides $v$ then considering the points of 
$\F{v/k}{q^k}$ as $k$-dimensional subspaces over $\mathbb{F}_q$ gives a construction of $k$-spreads. If $k$ does not divide $v$, then no 
$k$-spreads exist. Vector space partitions of type $k^{m_k}1^{m_1}$ are known under the name \emph{partial $k$-spreads}. More precisely, 
a partial $k$-spread in $\F{v}{q}$ is a set $\mathcal{K}$ of $k$-subspaces such that each point of the ambient space $\F{v}{q}$ is covered 
at most by one of its elements. Adding the set of uncovered points, which are also called \emph{holes}, gives a vector space partition of type 
$k^{m_k}1^{m_1}$. Maximizing $m_k=\# \mathcal{K}$ is equivalent to the minimization of $m_1$. If $d_1$ is the smallest dimension with 
$m_{d_1}\neq 0$, we call $m_{d_1}$ the \emph{length of the tail} and call the set of the corresponding $d_1$-subspace the \emph{tail}. Vector space 
partitions with a tail of small length are of special interest. In \cite{heden2009length} Olof Heden obtained:
\begin{Theorem} (Theorem 1 in \cite{heden2009length})
  \label{thm_length_of_tail}
  Let $\mathcal{P}$ be a vector space partition of type
  ${d_l}^{u_l}\ldots {d_2}^{u_2}{d_1}^{u_1}$ of $\SV{v}{q}$, where
  $u_1,u_2>0$ and $d_l>\dots>d_2>d_1\ge 1$.
  \begin{enumerate}[(i)]
  \item If $q^{d_2-d_1}$ does not divide $u_1$ and if $d_2<2d_1$, then
    $u_1\ge q^{d_1}+1$;
  \item if $q^{d_2-d_1}$ does not divide $u_1$ and if
    $d_2\ge 2d_1$, then either $d_1$ divides $d_2$ and
    $u_1=\gaussm{d_2}{1}{q}/\gaussm{d_1}{1}{q}$ 
    or $u_1>2q^{d_2-d_1}$;
  \item if $q^{d_2-d_1}$ divides $u_1$ and $d_2<2d_1$, then
    $u_1\ge q^{d_2}-q^{d_1}+q^{d_2-d_1}$;
  \item if $q^{d_2-d_1}$ divides $u_1$ and $d_2\ge 2d_1$, then
    $u_1\ge q^{d_2}$.
  \end{enumerate}   
\end{Theorem}
Moreover, in Theorem~2 and Theorem~3 he classified the possible sets of $d_1$-subspaces for $u_1=q^{d_1}+1$ and 
$u_1=\gaussm{d_2}{1}{q}/\gaussm{d_1}{1}{q}$, respectively. The results were obtained using the theory of mixed perfect $1$-codes, see 
e.g.\ \cite{herzog1972group}.  

In~\cite{blokhuis1989heden} the authors improved a lower bound of Heden on the size of inclusion-maximal partial $2$-spreads by translating the 
underlying techniques into geometry. Here we improve Theorem~\ref{thm_length_of_tail}(ii). The underlying geometric structure is the 
set $\mathcal{N}$ of $d_1$-subspaces of a vector space partition $\mathcal{P}$ of type ${d_l}^{u_l}\ldots {d_2}^{u_2}{d_1}^{u_1}$.  
For $d_1$ this is just a set of points in $\F{v}{q}$. It can be shown that the existence of $\mathcal{P}$ implies 
$\#\mathcal{N}\equiv \# \left(\mathcal{N}\cap H\right)\pmod {q^{d_2-1}}$ for every hyperplane $H$ of $\F{v}{q}$, see e.g.\ 
\cite{partial_spreads_and_vectorspace_partitions}. Taking a vector representation of the elements of $\mathcal{N}$ as columns of a
generator matrix, we obtain a corresponding (projective) linear code $\mathcal{C}$ over $\mathbb{F}_q$. The modulo constraints for 
$\mathcal{N}$ are equivalent to the property that the Hamming weights of the codewords of $\mathcal{C}$ are divisible by $q^{d_2-1}$. 
The study of so-called divisible codes, where the Hamming weights of the codewords of a linear code are divisible by some factor $\Delta>1$, was 
initiated by Harold Ward, see e.g.\ \cite{ward1981divisible}. The MacWilliams identities, linking the weight distribution of a linear code 
with the weight distribution of its dual code, can be relaxed to a linear program. Incorporating some information about the weight distribution 
of a linear code may result in an infeasible linear program, which then certifies the non-existence of such a code. This technique is known 
under the name linear programming method for codes and was more generally developed for association schemes by 
Philip Delsarte \cite{delsarte1973algebraic}. In~\cite{kurz2017packing} analytic solutions of linear programs for projective $q^r$-divisible 
linear codes have been applied in order to compute upper bounds for partial $k$-spreads. Indeed, all currently known upper bounds for 
partial $k$-spreads can be deduced from this method, see \cite{partial_spreads_and_vectorspace_partitions} for a survey.   

Here, we generalize the approach to the case $d_1>1$ by studying the properties of the set $\mathcal{N}$ of $d_1$-subspaces of 
a vector space partition $\mathcal{P}$ of $\F{v}{q}$ of type ${d_l}^{u_l}\ldots {d_2}^{u_2}{d_1}^{u_1}$ in Section~\ref{sec_divisible}. It 
turns out that we have $\#\mathcal{N}\equiv \#(\mathcal{N}\cap H)\pmod{q^{d_2-d_1}}$ for every hyperplane $H$ of $\SV{v}{q}$, see 
Lemma~\ref{lemma_divisible}, which we introduce as a definition of a \emph{$q^{d_2-d_1}$-divisible} set of $k$-subspaces with trivial 
intersection. By elementary counting techniques we obtain a partial substitute for the MacWilliams identities, see 
the equations (\ref{eq1}) and (\ref{eq2}). These imply some analytical criteria for the non-existence of such sets $\mathcal{N}$, which are 
used in Section~\ref{sec_reproof} to reprove Theorem~\ref{thm_length_of_tail}. By an improved analysis we tighten Theorem~\ref{thm_length_of_tail} 
to Theorem~\ref{main_thm}. More precisely, the second lower bound of Theorem~\ref{thm_length_of_tail}(ii) is improved. We close with some 
numerical results on the spectrum of the possible cardinalities of $\mathcal{N}$ and pose some open problems.

\section{Sets of disjoint $k$-subspaces and their incidences with hyperplanes}
\label{sec_divisible}
\noindent
For a positive integer $k$ let $\mathcal{N}$ be a set of pairwise disjoint, i.e., having trivial intersection, $k$-subspaces in $\F{v}{q}$, where we assume that the 
$k$-subspaces from $\mathcal{N}$ span $\F{v}{q}$, i.e., $v$ is minimally chosen. By $a_i$ we denote 
the number of hyperplanes $H$ of $\F{v}{q}$ with $\#(\mathcal{N}\cap H):=\#\{U\in\mathcal{N}\,:\,U\le H\}=i$ and set $n:=\#\mathcal{N}$. 
Due to our assumption on the minimality of the dimension $v$ not all $n$ elements from $\mathcal{N}$ can be contained in a hyperplane.  
Double-counting the incidences of 
the tuples $(H)$, $(B_1,H)$, and $(B_1,B_2,H)$, where $H$ is a hyperplane and $B_1\neq B_2$ are elements of $\mathcal{N}$ contained in $H$ gives:
\begin{equation}
  \sum_{i=0}^{n-1} a_i=\gaussm{v}{1}{q}, \quad
  \sum_{i=0}^{n-1} ia_i=n\cdot \gaussm{v-k}{1}{q}, \quad\text{and}\quad
  \sum_{i=0}^{n-1} i(i-1)a_i=n(n-1)\cdot \gaussm{v-2k}{1}{q}.\label{eq1}
\end{equation} 
For three different elements $B_1,B_2,B_3$ of $\mathcal{N}$ their span $\langle B_1,B_2,B_3\rangle$ has a dimension $i$ between $2k$ and $3k$. Denoting 
the number of corresponding triples by $b_i$, double-counting tuples $(B_1,B_2,B_3,H)$, where $H$ is a hyperplane and $B_1,B_2,B_3$ are pairwise different 
elements of $\mathcal{N}$ contained in $H$, gives:
\begin{equation}
  \sum_{i=0}^{n-1} i(i-1)(i-2)a_i=\sum_{i=2k}^{3k} b_i\gaussm{v-i}{1}{q} \quad\text{and}\quad \sum_{i=2k}^{3k} b_i =n(n-1)(n-2).\label{eq2}
\end{equation}
Given parameters $q$, $k$, $n$, and $v$ the so-called \emph{(integer) linear programming method} asks for a solution of the equation system 
given by (\ref{eq1}) and (\ref{eq2}) with 
$a_i,b_i\in\mathbb{R}_{\ge 0}$ ($a_i,b_i\in\mathbb{N}$). If no solution exists, then no corresponding set $\mathcal{N}$ can exist. For $k=1$ 
the equations from (\ref{eq1}) and (\ref{eq2}) correspond to the first four MacWilliams identities, see e.g.\ 
\cite{partial_spreads_and_vectorspace_partitions}.

If there is a single non-zero value $a_i$ the system can be solved analytically.
\begin{Lemma}
  \label{lemma_one_hyperplane_type}
  If $a_i=0$ for all $i\neq r>0$ and $k<v$ in the above setting, then there exists an integer $s\ge 2$ with $v=sk$ and $\mathcal{N}$  
  is a $k$-spread. Additionally we have $r=\frac{q^{v-k}-1}{q^k-1}$.
\end{Lemma}
\begin{Proof}
  Solving 
  (\ref{eq1}) for $r$, $a_r$, and $n$ gives $n=\frac{q^{2v-k}-q^{v}-q^{v-k}+1}{q^v-q^{v-k}-q^k+1}$.
  Writing $v=sk+t$ with $s,t\in\mathbb{N}$ and $0\le t<k$ we obtain
  $
    n=\sum_{i=1}^s q^{v-ik} +\frac{q^{v-k+t}-q^{v-k}-q^t+1}{q^v-q^{v-k}-q^k+1}
  $. Since $n\in\mathbb{N}$ and $0\le q^{v-k+t}-q^{v-k}-q^t+1<q^v-q^{v-k}-q^k+1$ we have $q^{v-k+t}-q^{v-k}-q^t+1=0$ so that $t=0$ 
  and $n=\frac{q^v-1}{q^k-1}$. Counting points gives that $\mathcal{N}$ partitions $\F{v}{q}$. 
\end{Proof}
We remark that $r=0$ forces $n\in\{0,1\}$ so that $\mathcal{N}$ is empty or consists of a single $k$-subspace in $\F{k}{q}$ and $v=k$ implies 
the latter case. So, these degenerated cases correspond to $s\in\{0,1\}$ in Lemma~\ref{lemma_one_hyperplane_type}. As pointed out after  
\cite[Theorem 2]{heden2009length}, such results can be proved in different ways. 
While the case that only one $a_i$ is non-zero is rather special, 
we can show that many $a_i$ are equal to zero in our setting.
\begin{Lemma}
  \label{lemma_divisible}
  Let $\mathcal{P}$ be a vector space partition of type ${d_l}^{u_l}\ldots {d_2}^{u_2}{d_1}^{u_1}$ of $\SV{v}{q}$, where
  $u_1,u_2>0$, and let $\mathcal{N}$ be the set of $d_1$-subspaces. Then, we have 
  $\#\mathcal{N}\equiv \#(\mathcal{N}\cap H)\pmod{q^{d_2-d_1}}$ for every hyperplane $H$ of $\SV{v}{q}$.
\end{Lemma}
\begin{Proof}
  For each $U\in\mathcal{P}$ we have $\dim(U\cap H)\in\{\dim(U),\dim(U)-1\}$. So counting points in $\SV{v}{q}$ and $H$ gives the existence 
  of integers $a,a'$ with $m\cdot\gaussm{d_2}{1}{q}+aq^{d_2}+u_1\gaussm{d_1}{1}{q}=\gaussm{v}{1}{q}$ and 
  $m\cdot\gaussm{d_2-1}{1}{q}+a'q^{d_2-1}+u_1'q^{d_1-1}+ u_1\gaussm{d_1-1}{1}{q}=\gaussm{v-1}{1}{q}$, where $m:=\sum_{i=2}^l u_i$ and 
  $u_1':=\#(\mathcal{N}\cap H)$. By subtraction we obtain $mq^{d_2-1}+aq^{d_2}-a'q^{d_2-1}+u_1q^{d_1-1}-u_1'q^{d_1-1}=q^{v-1}$, so that 
  $u_1q^{d_1-1}\equiv u_1'q^{d_1-1}\pmod {q^{d_2-1}}$.
\end{Proof}

\begin{Definition}
  Let $\mathcal{N}$ be a set of $k$-subspaces in $\SV{v}{q}$. 
  If there exists 
  a positive integer $r$ such that $a_i$ is non-zero only if $\#\mathcal{N}-i$ is divisible by $q^{r}$ and the $k$-subspaces are pairwise disjoint, 
  then we call $\mathcal{N}$ \emph{$q^r$-divisible}.  
\end{Definition}

Using the notation of Lemma~\ref{lemma_divisible}, $\mathcal{N}$ is $q^{d_2-d_1}$-divisible. As mentioned in the introduction, for $d_1=1$, 
taking the elements of $\mathcal{N}$ as columns of a generator matrix, we obtain a projective linear code, whose Hamming weights are 
divisible by $q^{d_2-1}$.

\begin{Example}
  \label{construction_1}
  For integers $k\ge 2$ and $r=ak+b$ with $0\le b<k$ let $\mathcal{N}$ be a $k$-spread of $\F{(a+2)k}{q}$. Starting from 
  a $(a+2)k$-spread in $\F{2(a+2)k}{q}$ we obtain a vector space partition $\mathcal{P}$ by replacing one $(a+2)k$-dimensional 
  spread element with $\mathcal{N}$. From Lemma~\ref{lemma_divisible} and $q^r | q^{(a+2)k-k}=q^{(a+1)k}$ we deduce that the set 
  $\mathcal{N}$ of $k$-subspaces is $q^r$-divisible. Its cardinality is given by $\gaussm{(a+2)k}{1}{q}/\gaussm{k}{1}{q}$.
\end{Example}

\begin{Example}
  \label{construction_2}
  For integers $k\ge 2$ and $r\ge 1$ let $n=k+r$ and consider a matrix representation
  $M\colon\mathbb{F}_{q^n}\to\mathbb{F}_q^{n\times n}$ of $\mathbb{F}_{q^n}/\mathbb{F}_q$, obtained by
  expressing the multiplication maps $\mu_\alpha\colon\mathbb{F}_{q^n}\to\mathbb{F}_{q^n}$, $x\mapsto\alpha x$, 
  which are linear over $\mathbb{F}_q$, in terms of a fixed basis of $\mathbb{F}_{q^n}/\mathbb{F}_q$. Then, all
  matrices in $M(\mathbb{F}_{q^n})$ are invertible and have mutual rank distance 
  $d_{\operatorname{R}}(A,B):=\operatorname{rk}(A-B)=n$, see e.g.\ \cite{partial_spreads_and_vectorspace_partitions} 
  for proofs of these and the subsequent facts. In other words, the matrices of $M(\mathbb{F}_{q^n})$ form a maximum rank distance 
  code with minimum rank distance $n$ and cardinality $q^n$. 
  
  Now let $\mathcal{B}\subseteq\mathbb{F}_q^{k\times n}$ be the matrix code
  obtained from $M(\mathbb{F}_{q^n})$ by deleting the last $n-k$ rows, say, of
  every matrix. Then $\mathcal{B}$ has cardinality minimum rank distance $k$.  
  Hence, by applying the lifting construction $B\mapsto(I_k|B)$, where $I_k$ is the $k\times k$ identity matrix,  
  to $\mathcal{B}$ we obtain a partial $k$-spread $\mathcal{N}$ in $\mathbb{F}_q^v$ of size $q^n=q^{k+r}$. 
  Since precisely the points outside the ($k+r$)-subspace $S=\left\{x\in\mathbb{F}_q^v\,:\,x_1=x_2=\dots=x_k=0\right\}$ 
  are covered, $\mathcal{P}=\mathcal{N}\cup\{S\}$ is a vector space partition of $\F{2k+r}{q}$ and $\mathcal{N}$ 
  is $q^{k+r}$-divisible with cardinality $q^{k+r}$. 
\end{Example}

From the first two equations of~(\ref{eq1}) we deduce:
\begin{Lemma}
  \label{lemma_linear_condition}
  For a $q^r$-divisible set $\mathcal{N}$ of $k$-subspaces in $\SV{v}{q}$, there exists a hyperplane $H$ with 
  $\#(\mathcal{N}\cap H)\le n/q^k$.
\end{Lemma}
\begin{Proof}
  Let $i$ be the smallest index with $a_i\neq 0$. 
  Then, the first two equations of~(\ref{eq1}) are equivalent to 
  $\sum_{j\ge 0} a_{i+q^rj}=\gaussm{v}{1}{q}$ and $\sum_{j\ge 0} \left(i+q^rj\right)\cdot a_{i+q^rj}=n\gaussm{v-k}{1}{q}$. 
  Subtracting 
  $i$ times the first equation from the second equation gives $\sum_{j>0} q^rja_{i+q^rj}=n\cdot\frac{q^{v-k}-1}{q-1}-i\cdot\frac{q^v-1}{q-1}$. 
  Since the left-hand side is non-negative, we have $i\le \frac{q^{v-k}-1}{q^v-1}\cdot n\le \frac{n}{q^k}$. 
\end{Proof}

Stated less technical, the proof of Lemma~\ref{lemma_linear_condition} is given by the fact that the hyperplane with the minimum number of $k$-subspaces 
contains at most as many $k$-subspaces as the average number of $k$-subspaces per hyperplane.

Taking also the third equation of~(\ref{eq1}) into account implies a quadratic criterion:
\begin{Lemma}
  \label{lemma_quadratic_condition}
  Let $m\in \mathbb{Z}$ and $\mathcal{N}$ be a $q^r$-divisible set of $k$-subspaces in $\SV{v}{q}$. Then, 
  $\tau(n,q^r,q^k,m)\cdot q^{v-2k-2r}-m(m-1)\ge 0$, where $\tau(n,\Delta,u,m):=\Delta^2u^2m(m-1)-n(2m-1)u(u-1)\Delta+n(u-1)(n(u-1)+1)$. 
\end{Lemma}
\begin{Proof}
  With $y=q^{v-2k}$, $u=q^k$, and $\Delta=q^r$, we can rewrite the equations of~(\ref{eq1}) to
  $u^2y-1=(q-1)\sum_{i\in\mathbb{Z}} a_i$, $n\cdot\left(uy-1\right)=(q-1)\sum_{i\in\mathbb{Z}}ia_i$, and
  $n(n-1)\cdot\left(y-1\right)=\sum_{i\in\mathbb{Z}} i(i-1)a_i$. $(n-m\Delta)(n-(m-1)\Delta)$ times the first minus 
  $2n-(2m-1)\Delta-1$ times the second plus the third equation gives $y\cdot \tau(n,\Delta,u,m)-\Delta^2m(m-1)
  =(q-1)\sum_{i\in\mathbb{Z}} (n-m\Delta-i)(n-(m-1)\Delta-i)a_i=(q-1)\sum_{h\in\mathbb{Z}} \Delta^2(m-h)(m-h+1)a_{n-h\Delta}\ge 0$.  
\end{Proof}

As a preparation we present another classification result: 
\begin{Lemma}
  \label{lemma_class_2}
   If $\mathcal{N}$ is a $q$-divisible set of $k$-subspaces in $\SV{v}{q}$ of cardinality $q^k+1$, then $\mathcal{N}$ partitions $\SV{2k}{q}$.  
\end{Lemma}
\begin{Proof}
  Setting $c_i:=(q-1)a_{1+iq}$ and $l:=q^{k-1}-1$ we can rewrite the equations of~(\ref{eq1}) to
  $\sum_{i=0}^{l} c_i=q^v-1$,
  $\sum_{i=0}^{l} (1+iq)c_i=(q^k+1)\left(q^{v-k}-1\right)$, and 
  $\sum_{i=0}^{l} iq(1+iq)c_i=(q^k+1)q^k\left(q^{v-2k}-1\right)$. Since $ql+1$ times the second minus $ql+1$ times the first minus the 
  third equation gives $0\le \sum_{i=0}^l iq^2(l-i)c_i=-q^{k+1}\left(q^{v-2k}-1\right)$, we have $v=2k$. 
  Every point of $\SV{v}{q}$ is covered by an element from $\mathcal{N}$ due to $\gaussm{2k}{1}{q}/\gaussm{k}{1}{q}=q^k+1$.  
\end{Proof}

\section{Proof of Heden's results and further improvements}
\label{sec_reproof}

\noindent
Let $\mathcal{P}$ be a vector space partition of type ${d_l}^{u_l}\ldots {d_2}^{u_2}{d_1}^{u_1}$ of $\SV{v'}{q}$, where
$u_1,u_2>0$, $d_l>\dots>d_2>d_1\ge 1$. Let $\mathcal{N}$ be the set of $d_1$-subspaces and $V$ be the subspace spanned by $\mathcal{N}$. 
By $n$ we denote the cardinality of $\mathcal{N}$ and by $a_i$ we denote the number of hyperplanes of $V$ that contain exactly 
$i$ elements from $\mathcal{N}$. 

Assume that $q^{d_2-d_1}$ does not divide $u_1$. We have $\#(\mathcal{N}\cap H)\ge 1$ for every hyperplane $H$ of $V$ due 
to Lemma~\ref{lemma_divisible}, so that Lemma~\ref{lemma_linear_condition} gives $u_1\ge q^{d_1}$. Thus, we have $u_1\ge q^{d_1}+1$. 
If $u=q^{d_1}+1$ then we can apply Lemma~\ref{lemma_class_2} for the classification of the possible sets $\mathcal{N}$. If 
$u_1<2q^{d_2-d_1}$ then for $a_i>0$ we have $i<q^{d_2-d_1}$ and $i\equiv u_1\pmod {q^{d_2-d_1}}$ so that we can apply 
Lemma~\ref{lemma_one_hyperplane_type}. Thus, either $d_2$ divides $d_1$ and $u_1=(q^{d_2}-1)/(q^{d_1}-1)$ or $u_1>2q^{d_2-d_1}$. 
The first case can be attained by a $d_2$-spread where one $d_2$-subspace is replaced by a $d_1$-spread, see Example~\ref{construction_1}. 
We remark that no assumption on the relation between $d_2$ and $d_1$ is used in our derivation. However, if $d_2<2d_1$ then $d_1$ cannot 
divide $d_2$ and $q^d_1+1>2q^{d_2-d_1}$.  

Assume that $q^{d_2-d_1}$ divides $u_1$. Setting $\Delta=q^{d_2-d_1}$, $u=q^{d_1}$, $n=\Delta l$, and $m=l$\footnote[2]{The choice for $m$ can be 
obtained by minimizing $\tau(n,\Delta,u,m)$, i.e., solving $\frac{\partial \tau(n,\Delta,u,m)}{\partial m}=0$ and rounding.} for some 
integer $l$, we conclude $\tau(n,\Delta,u,m)=\Delta l(\Delta l-\Delta u +u-1)\ge 0$ from Lemma~\ref{lemma_quadratic_condition}, 
so that $l\ge \left\lceil u-\frac{u}{\Delta}+\frac{1}{\Delta}\right\rceil$. The right-hand side is equal to $u=q^{d_1}$ if $d_2\ge 2d_1$ and 
to $u-u/\Delta+1=q^{d_1}-q^{2d_1-d_2}+1$ otherwise, which is equivalent to $n\ge q^{d_2}$ and $n\ge q^{d_2}-q^{d_1}+q^{d_2-d_1}$. We remark 
that equality is achievable in the latter case via the $2$-weight codes constructed in \cite{bierbrauer1997family} (with parameters 
$n'=d_1$ and $m=d_2-d_1$). We do not know whether the corresponding $q^{d_2-d_1}$-divisible set of $d_1$-subspaces can be realized as 
a vector space partition of $\SV{v}{q}$.\footnote[3]{A suitable test case might be to decide whether a vector space partition of type 
$4^4 3^{135} 2^6$ exists in $\SV{10}{2}$.} For the first case see Example~\ref{construction_2}.

The above comprises \cite[Theorems 1-4]{heden2009length}. Given the stated examples, just Theorem~\ref{thm_length_of_tail}(ii), for the case 
where $d_1$ does not divide $d_2$, leaves some space for improving the lower bound on $u_1$. To that end we analyze 
Lemma~\ref{lemma_quadratic_condition} in more detail. Since the statements look rather technical and complicated we first give 
a justification for the necessity of this fact. Via the quadratic inequality of Lemma~\ref{lemma_quadratic_condition} intervals 
of cardinalities can be excluded for different values of the parameter $m$. However, some cardinalities are indeed feasible. If $r=ak+b$ 
with $0\le b<k$ then the two constructions from Example~\ref{construction_1} and Example~\ref{construction_2} give $q^r$-divisible set 
of $k$-subspaces of cardinality $\gaussm{(a+2)k}{1}{q}/\gaussm{k}{1}{q}$ and $q^{k+r}$, respectively. For $q=2$, $r=3$, $k=2$ the cardinalities 
of these two examples are given by $21$ and $32$. In general, each two $q^r$-divisible sets $\mathcal{N}_1$ and $\mathcal{N}_2$ of $k$-subspaces can be 
combined to a $q^r$-divisible set of $k$-subspaces of cardinality $\#\mathcal{N}_1+\#\mathcal{N}_2$. Since 
$\gaussm{(a+2)k}{1}{q}/\gaussm{k}{1}{q}$ and $q^{k+r}$ are coprime there exists some integer $F_q(k,r)$ such that $q^r$-divisible sets 
of $k$-subspaces exist for every cardinality $n>F_q(k,r)$. Below that number some cardinalities can be excluded, but their density decreases 
with increasing $n$. Our numerical example is continued after the proof of Theorem~\ref{main_thm}.

\begin{Proposition}
  \label{prop_analysis_quadratic}
  Let $\mathcal{N}$ be a $q^r$-divisible set of $k$-subspaces in $\SV{v}{q}$, $u=q^k$ and $\Delta=q^r$. Then, 
  $n\notin \Big[1,\frac{q^{k+r}-1}{q^r-1}\Big)$ and $$n\notin\left[
  \left\lceil\frac{1}{u-1}\cdot\left(\Delta um -\frac{\Delta u+1}{2}-\frac{1}{2}\sqrt{\omega}\right)\right\rceil,
  \left\lfloor\frac{1}{u-1}\cdot\left(\Delta um -\frac{\Delta u+1}{2}+\frac{1}{2}\sqrt{\omega}\right)\right\rfloor\right],$$ 
  where $\omega=\left(\Delta u-2m\right)^2+\left(2\Delta u+1-4m^2\right)$, for all $m\in\mathbb{N}$ with 
  $2\le m\le \left\lfloor\frac{\Delta u}{4}+\frac{1}{2}+\frac{1}{4\Delta u} \right\rfloor$.
\end{Proposition}
\begin{Proof}
  We set $\overline{\Delta}=\Delta u$ and $\overline{n}=n(u-1)$ so that 
  $\tau(n,\Delta,u,m)=\overline{\Delta}^2m(m-1)-\overline{n}\overline{\Delta}(2m-1)+\overline{n}(\overline{n}+1)$. We have 
  $\tau(n,\Delta,u,m)\le 0$ iff 
    $\left|\overline{n}-\overline{\Delta}m+\frac{\overline{\Delta}+1}{2} \right|\le \frac{1}{2}\sqrt{\overline{\Delta}^2-4m\overline{\Delta}+2\overline{\Delta}+1}$
  and $m\le \frac{\overline{\Delta}}{4}+\frac{1}{2}+\frac{1}{4\overline{\Delta}}$. Rewriting and applying Lemma~\ref{lemma_quadratic_condition} 
  with $1\le m\le \left\lfloor\frac{\Delta u}{4}+\frac{1}{2}+\frac{1}{4\Delta u} \right\rfloor$ gives the result since $m(m-1)>0$ for $m\ge 2$.
\end{Proof}

\begin{Proposition}
  Let $\mathcal{N}$ be a $q^r$-divisible set of $k$-subspaces in $\SV{v}{q}$, where $r=ak+b$ with $a,b\in\mathbb{N}$, $0<b<k$ and $a\ge 1$. Then, 
  $n\ge \frac{q^{(a+2)k}-1}{q^k-1}=q^r\cdot q^{k-b}+\frac{q^r\cdot q^{k-b}-1}{q^k-1}=\Delta q^{k-b}+q^k\Theta+1$, where $\Delta:=q^r$ and 
  $\Theta:=\frac{q^{ak}-1}{q^k-1}$.  
\end{Proposition}
\begin{Proof}
  From Lemma~\ref{lemma_one_hyperplane_type} we conclude $n\ge 2q^r$ and set $u=q^k$. For $2\le m\le q^{k-b}$ we have 
  $2\Delta u+1-4m^2>0$, so that Proposition~\ref{prop_analysis_quadratic} gives $n\notin \left[
  \left\lceil\frac{\Delta u(m-1)-1/2+m}{u-1}\right\rceil,
  \left\lfloor\frac{\Delta um-1/2-m}{u-1}\right\rfloor
  \right]$. Since $\Delta(m-1)\le \left\lceil\frac{\Delta u(m-1)-1/2+m}{u-1}\right\rceil
  =\Delta (m-1)+\left\lceil\frac{\Delta (m-1)-1/2+m}{u-1}\right\rceil\le \Delta m$ and 
  $\left\lfloor\frac{\Delta um-1/2-m}{u-1}\right\rfloor=\Delta m+ mq^b\Theta+\left\lfloor\frac{mq^b-1/2-m}{q^k-1}\right\rfloor=\Delta m+ mq^b\Theta$, 
  we conclude $n\notin\left[\Delta m,\Delta m+ mq^b\Theta\right]$ for $2\le m\le q^{k-b}$.
  
  It remains to show $n\notin\left[\Delta m,\Delta m+ mq^b\Theta+1,\Delta(m+1)-1\right]=:I_m$ for all $2\le m\le q^{k-b}-1$. If $n\in I_m$, then
  we can write $n=\Delta m+mq^b\Theta+x$ with $x\ge 1$ and $mq^b\Theta+x<\Delta$, so that $q^k\cdot \left(mq^b\Theta+x\right) =
  \Delta m + mq^b\Theta +\left(xq^k-mq^b\right)<\Delta m + mq^b\Theta+x=n$, which contradicts Lemma~\ref{lemma_linear_condition}.  
\end{Proof}

In other words, in the case of Theorem~\ref{thm_length_of_tail}(ii), where $d_2=ad_1+b$ with $0<b<d_1$ and $a,b\in\mathbb{N}$, we have 
$u_1\ge q^{d_2-d_1}\cdot q^{d_1-b}+\frac{q^{(a+1)d_1}-1}{q^{d_1}-1}=\frac{q^{(a+2)d_1}-1}{q^{d_1}-1}$, which can be attained by an $d_1$-spread 
in $\SV{(a+2)d_1}{q}$. Without the knowledge of $b$, we can state $u_1\ge q\cdot q^{d_2-d_1}+\left\lceil\frac{q^{d_2+1}-1}{q^{d_1}-1}\right\rceil$, 
which also improves Theorem~\ref{thm_length_of_tail}(ii) and is tight whenever $d_2+1$ is divisible by $d_1$. Summarizing our findings we obtain our 
main theorem:
\begin{Theorem}
  \label{main_thm}
  For a non-empty $q^r$-divisible set $\mathcal{N}$ of $k$-subspaces in $\SV{v}{q}$ 
  the following bounds on $n=\#\mathcal{N}$ are tight.
  \begin{enumerate}[(i)]
    \item We have $n\ge q^k+1$ and if $r\ge k$ then either $k$ divides $r$ and $n\ge \frac{q^{k+r}-1}{q^k-1}$ or $n\ge\frac{q^{(a+2)k}-1}{q^k-1}$, 
          where $r=ak+b$ with $0<b<k$ and $a,b\in\mathbb{N}$.
    \item Let $q^r$ divide $n$. If $r<k$ then $n\ge q^{k+r}-q^k+q^r$ and $n\ge q^{k+r}$ otherwise. 
  \end{enumerate}  
\end{Theorem} 
For (i) the lower bounds are attained by $k$-spreads, see Example~\ref{construction_1}. For (ii) the second lower bound is attained 
by a construction based on lifted MRD codes, see Example~\ref{construction_2}. In the other case the $2$-weight codes constructed in 
\cite{bierbrauer1997family} attain the lower bound. Thus, Theorem~\ref{main_thm} is tight and implies an improvement of 
Theorem~\ref{thm_length_of_tail}(ii).   

While the smallest cardinality of a non-empty $q^r$-divisible set of $k$-subspaces over $\mathbb{F}_q$ has been determined, 
the spectrum of possible cardinalities remains widely unknown. 
For $k=1$ \cite[Theorem 12]{partial_spreads_and_vectorspace_partitions} states that 
either $n>rq^{r+1}$ or there exist integers $a,b$ with $n=a\gaussm{r+1}{1}{q}+bq^{r+1}$ and bounds for the maximum excluded cardinality 
have been determined in \cite{projective_divisible_binary_codes}. However, Lemma~\ref{lemma_linear_condition} and Lemma~\ref{lemma_quadratic_condition}, applied 
via Proposition~\ref{prop_analysis_quadratic}, give  restrictions going far beyond Theorem~\ref{main_thm}. For $q=2$, $r=3$, $k=2$, and $n\le 81$ we 
exemplarily state that only $n\in\{21,31,32,33,42,$ $43,44,52,\ldots,55,62,\ldots,66,72,\ldots,78\}$ might be attainable. 
The mentioned constructions cover 
the cases $n\in\{21,32,42,53,63,64,74\}\subseteq \{21a+32b\,:\, a,b\in\mathbb{N}\}$. Replacing the lines by their contained $3$ points, we obtain 
$2^4$-divisible sets of $1$-subspaces in $\SV{v}{q}$ of cardinality $3n$, for which two further exclusion criteria have been presented in 
\cite{partial_spreads_and_vectorspace_partitions}, excluding the cases $n\in\{33,44\}$. 
\cite[Lemma~23]{partial_spreads_and_vectorspace_partitions} is based on a cubic polynomial obtained from (\ref{eq1}) and (\ref{eq2}), similar to 
the quadratic polynomial from Lemma~\ref{lemma_quadratic_condition} obtained from (\ref{eq1}). Here, the presence of $k$ additional $b_i$-variables 
may make the analysis more difficult for $k>1$. For a $q^r$-divisible set $\mathcal{N}$ of $1$-subspaces we have that $\mathcal{N}\cap H$ is 
$q^{r-1}$-divisible for every hyperplane $H$, which allows a recursive application of the linear programming method. For $k>1$ we need to 
consider $k$-subspaces and $k-1$-subspaces in $H$, see \cite[Section 6.3]{partial_spreads_and_vectorspace_partitions}, which makes the bookkeeping more 
complicated. 

The determination of the possible spectrum of cardinalities of $q^r$-divisible sets of $k$-subspaces remains an interesting open problem. Even for 
small parameters this might be challenging. A possible intermediate step is the determination of the number $F_q(k,r)$ being similar to the Frobenius 
number. Extending the small list of constructions is also worthwhile.  

\section*{Acknowledgement}
I am very thankful for the comments of two anonymous reviewers, which helped to improve the paper. 


\end{document}